\documentclass[12pt]{article}
\usepackage{amsmath}
\usepackage{amsmath,amssymb, amsthm, latexsym, amsfonts, epsfig, color,graphicx,}
\usepackage{mathrsfs}
\usepackage{enumitem}
\usepackage{indentfirst}
\usepackage{chngcntr}
\usepackage{bbm}
\usepackage[blocks]{authblk} 
\usepackage{newtxmath} 
\usepackage[colorlinks=True,linkcolor=blue,anchorcolor=blue,citecolor=blue,CJKbookmarks=true]{hyperref}
\usepackage{geometry}
\usepackage{verbatim}
\usepackage{amsfonts}
\DeclareMathOperator{\supp}{supp}

\begin{document}
\renewcommand{\a}{\alpha}
\newcommand{\D}{\Delta}
\newcommand{\ddt}{\frac{d}{dt}}
\newcommand{\R}{\mathbb R}
\newcommand{\C}{\mathbb C}
\newcommand{\Sph}{\mathbb S^{d-1}}
\newcommand{\Lip}{\operatorname{Lip}}
\newcommand{\BM}{\mathrm{BM}}
\newcommand{\dd}{\,\mathrm d}
\newcommand{\la}{\langle}
\newcommand{\ra}{\rangle}
\counterwithin{equation}{section}
\newcommand{\e}{\epsilon}
\newcommand{\eps}{\varepsilon}
\newtheorem{theorem}{Theorem}[section]
\newtheorem{proposition}{Proposition}[section]
\newtheorem{lemma}[theorem]{Lemma}
\newtheorem{remark}[theorem]{Remark}
\newtheorem{example}{Example}[section]
\newtheorem{definition}{Definition}[section]
\newtheorem{corollary}[theorem]{Corollary}
\newtheorem{citedtheorem}{Theorem}
\renewcommand{\thecitedtheorem}{\Alph{citedtheorem}}
\makeatletter
\newcommand{\rmnum}[1]{\romannumeral #1}
\newcommand{\Rmnum}[1]{\expandafter\@slowromancap\romannumeral #1@}
\makeatother

\title{Radial Beurling--Malliavin multipliers:\\
a characterization and sharp logarithmic criteria}
\author{Longben Wei
}
\date{}
\maketitle
\begin{abstract} 
We give a complete characterization of the possibly nonradial weights that admit radial Beurling--Malliavin multipliers in $\R^d$, $d\ge2$. For a nonnegative globally Lipschitz function $\Omega$, such multipliers exist at every positive bandwidth if and only if
\[
 \int_0^\infty
 \frac{\max_{|x|=r}\Omega(x)}{1+r^2}\dd r<\infty.
\]
Thus a single spherical maximum condition exactly determines
whether the prescribed bound $|f|\le e^{-\Omega}$ can be attained
by nonzero radial bandlimited $L^2$ functions.

Building on this characterization, we establish two sharp
logarithmic results. First, we replace the power correction in Vasilyev's
volume-integral criterion by an optimal logarithmic correction: finiteness of
$\int_{\R^d}\Omega(x)[\log(e+|x|)]^\beta
(1+|x|)^{-d-1}\dd x$
guarantees radial BM multipliers whenever $\beta>d-1$. For every $\beta\le d-1$, smooth nonnegative globally Lipschitz counterexamples satisfy this integral condition but admit no nonzero bandlimited multiplier, even without radiality. Second, under Cohen's directional integrability condition and Lipschitz regularity alone, radial multipliers exist at every positive bandwidth with
$|f(x)|\le\exp(-\Omega(x)/[\log(e+|x|)]^\gamma)$
for every $\gamma>1/2$. The threshold $1/2$ is optimal within the class of radial multipliers. Finally, we give a Hankel--Paley--Wiener proof of Vasilyev's radial theorem that treats all dimensions uniformly.

\end{abstract}
\section{Introduction}

The Beurling--Malliavin multiplier problem concerns the existence of nonzero functions with prescribed pointwise decay and arbitrarily small Fourier support. In several dimensions, a natural question is which weights admit multipliers that are radial, even when the weights themselves are not. We give a complete answer for weights $e^{-\Omega}$ with $\Omega\ge0$ globally Lipschitz: radial multipliers exist at every positive bandwidth if and only if the spherical maximum of $\Omega$ is integrable against
$(1+r^2)^{-1}\dd r$ on $[0,\infty)$. This characterization leads to sharp logarithmic refinements of Vasilyev's weighted volume criterion \cite{Vas} and an optimal logarithmic loss under Cohen's directional integrability condition \cite{Cohen}.

Throughout the paper, we use the Fourier transform
\[
 \widehat f(\xi)
 =\int_{\R^d}f(x)e^{-2\pi ix\cdot\xi}\dd x,
\]
extended to $L^2(\R^d)$ by Plancherel's theorem,
with the same normalization as in \cite{Vas}.
\begin{definition}
A weight $\omega:\R^d\to(0,1]$ is called a BM majorant if, for every $\sigma>0$, there exists $0\ne f\in L^2(\R^d)$ such that
\begin{equation}\label{eq:bm}
 \supp\widehat f\subset\overline B(0,\sigma),
 \qquad |f(x)|\le\omega(x)\quad(x\in\R^d).
\end{equation}
We call such a function $f$ a BM multiplier for $\omega$. The pointwise assertion refers to the continuous representative of $f$. We say that $\omega$ admits radial BM multipliers if $f$ can be chosen radial for every $\sigma>0$. The weight $\omega$ itself may be nonradial.
\end{definition}

\begin{samepage}
\begin{citedtheorem}[Classical Beurling--Malliavin multiplier theorem {\cite{BM}}]
\label{thm:classical-bm}
Let $\omega:\R\to(0,1]$ and write $\Omega=\log(1/\omega)$. If $\Omega$ is globally Lipschitz and
\[
 \int_{\R}\frac{\Omega(t)}{1+t^2}\dd t<\infty,
\]
then $\omega$ is a BM majorant on $\R$.
\end{citedtheorem}
\end{samepage}

The Beurling--Malliavin multiplier theorem is a cornerstone of one-dimensional harmonic analysis, revealing a deep connection between pointwise decay and spectral localization. For background and further developments, we refer to \cite[Part~II]{HavinJoricke} and the survey
\cite{MashreghiNazarovHavin}. A multiplier furnished by this theorem can impose decay adapted to a given set while increasing the Fourier radius by an arbitrarily small amount. Bourgain and Dyatlov use this construction in their proof of a fractal uncertainty principle and the resulting essential spectral gap for convex co-compact hyperbolic surfaces \cite{BD}. Han and Schlag extend this approach to a higher-dimensional fractal uncertainty principle with product structure on the Fourier side \cite{HanSchlag}. For further developments concerning fractal uncertainty principles and quantitative uniqueness, see \cite{JinZhang,JayeMitkovski}; for a generalization of the one-dimensional BM majorant theorem, see \cite{VasMajorant}.

The search for a several-variable multiplier theorem goes back to H\"ormander's work on convolution equations in convex domains \cite{Hormander}. Cohen establishes a higher-dimensional fractal uncertainty principle for sets that are porous on balls in physical space and on lines in frequency space. The key ingredient is a multidimensional Beurling--Malliavin multiplier theorem under directional integrability and symbol regularity hypotheses \cite{Cohen}. Vasilyev proves the radial theorem and derives a nonradial result under a single weighted volume condition \cite{Vas}. Our first contribution is to reformulate Vasilyev's radial theorem as a complete characterization of the general, possibly nonradial weights that admit radial BM multipliers. We then establish a sharp logarithmic refinement of his nonradial volume criterion. A further contribution is a radial multiplier theorem under Cohen's directional integrability condition and Lipschitz regularity alone, with an optimal logarithmic loss. We also give a Hankel--Paley--Wiener proof of Vasilyev's radial theorem that treats all dimensions uniformly. We begin with the radial theorem.

\begin{citedtheorem}[Vasilyev's radial theorem {\cite[Theorem 1]{Vas}}]
\label{thm:vas-radial}
Let $d\ge2$ and let $\Omega:\R^d\to[0,\infty)$ be radial and globally Lipschitz. If
\[
 \int_{\R^d}\frac{\Omega(x)}{(1+|x|^2)^{(d+1)/2}}\dd x<\infty,
\]
then $e^{-\Omega}$ admits radial BM multipliers.
\end{citedtheorem}

Vasilyev's radial theorem provides a sufficient condition for the existence of radial multipliers when the majorant is radial. A natural question is whether one can obtain a necessary and sufficient condition without assuming radiality of the majorant. As we shall see, such a characterization follows by considering the spherical maximum of the logarithmic weight. For a nonnegative globally Lipschitz function $\Omega$, define
\begin{equation}\label{eq:envelope-def}
 m_\Omega(r)=\max_{|x|=r}\Omega(x),\qquad
 \mathcal J(\Omega)=\int_0^\infty\frac{m_\Omega(r)}{1+r^2}\dd r.
\end{equation}
The following theorem gives the complete criterion for the existence of radial multipliers under a possibly nonradial weight.

\begin{theorem}\label{thm:main}
Let $d\ge2$ and let $\Omega:\R^d\to[0,\infty)$ be globally Lipschitz. Neither $\Omega$ nor the weight $\omega=e^{-\Omega}$ is assumed radial. The following conditions are equivalent.
\begin{enumerate}[label=\textup{(\roman*)},leftmargin=*]
\item $\mathcal J(\Omega)<\infty$.
\item There exist $\sigma_0>0$ and a nonzero radial $f\in L^2(\R^d)$ such that $\supp\widehat f\subset\overline B(0,\sigma_0)$ and $|f(x)|\le e^{-\Omega(x)}$ for all $x\in\R^d$.
\item For every $\sigma>0$, there is a nonzero radial $f_\sigma\in L^2(\R^d)$ such that $\supp\widehat f_\sigma\subset\overline B(0,\sigma)$ and $|f_\sigma(x)|\le e^{-\Omega(x)}$ for all $x\in\R^d$.
\end{enumerate}
\end{theorem}
Theorem~\ref{thm:main} characterizes the existence of radial multipliers satisfying $|f|\le e^{-\Omega}$ at every positive bandwidth. For weights with $\mathcal J(\Omega)=\infty$, any nonzero bandlimited multiplier, if one exists, must be nonradial. The proof of Theorem~\ref{thm:main} is quite simple. For sufficiency, the spherical maximum $m_\Omega(|x|)$ remains globally Lipschitz, so Theorem~\ref{thm:vas-radial} applies.
For necessity, a radial multiplier is bounded by $e^{-m_\Omega(r)}$ on each sphere. Bergman's necessary condition
for radial weights \cite[Proposition~1]{Bergman} then gives $\mathcal J(\Omega)<\infty$. 
\begin{remark}\label{rem:main}
Our integrability condition~\eqref{eq:envelope-def} is close in form to Cohen's. A comparison with his multiplier theorem \cite[Theorem~1.4]{Cohen} therefore clarifies the scope of Theorem~\ref{thm:main}. Cohen's directional growth condition is
\begin{equation}\label{eq:cohen-directional}
 G_\Omega(x)=\int_{1/2}^{2}\Omega(sx)\dd s,\qquad
 G_\Omega^*(r)=\max_{|x|=r}G_\Omega(x),\qquad
 \mathcal C_{\rm dir}(\Omega)
 =\int_0^\infty\frac{G_\Omega^*(r)}{1+r^2}\dd r<\infty.
\end{equation}
For radial majorants, this condition is equivalent to $\mathcal J(\Omega)<\infty$. For general nonradial majorants, however,
\[
 G_\Omega^*(r)
 =\sup_{\theta\in\Sph}\int_{1/2}^{2}\Omega(sr\theta)\dd s
 \le\int_{1/2}^{2}m_\Omega(sr)\dd s,
\]
and consequently
\[
 \mathcal C_{\rm dir}(\Omega)\le3\mathcal J(\Omega).
\]
The converse fails even for smooth weights satisfying Cohen's symbol estimates through order three; see
Section~\ref{sec:cohen-comparison}. Thus Cohen's integrability condition is strictly weaker. For weights with $\mathcal J(\Omega)=\infty$, any multiplier supplied by his theorem must be nonradial.

Although Cohen's integrability condition is weaker, his theorem imposes stronger regularity assumptions, whereas
Theorem~\ref{thm:main} requires only Lipschitz continuity. Whether Cohen's theorem remains valid under Lipschitz regularity alone remains an open question.
\end{remark}

Although Theorem~\ref{thm:main} already applies to general nonradial weights, its criterion involves the maximum of $\Omega$ on each sphere. A criterion expressed through a single weighted volume integral can be
easier to verify, especially when integral estimates for $\Omega$ are available. Vasilyev derives such a sufficient condition from his radial theorem \cite[Theorem~2]{Vas}. Our next goal is to sharpen this volume criterion by replacing the additional power in its weight with an optimal logarithmic correction. In our notation, Vasilyev's result reads as follows.

\begin{citedtheorem}[Vasilyev's nonradial criterion {\cite[Theorem 2]{Vas}}]
\label{thm:vas-nonradial}
Let $d\ge2$ and let $\Omega:\R^d\to[0,\infty)$ be globally Lipschitz. If
\begin{equation}\label{eq:old-power}
 \int_{\R^d}\frac{\Omega(x)}{(1+|x|)^{d+1-\varepsilon}}\dd x<\infty
 \quad\text{for some }\varepsilon>0,
\end{equation}
then $e^{-\Omega}$ admits radial BM multipliers.
\end{citedtheorem}

Our search for a logarithmic refinement is motivated by Bourgain and Dyatlov's proof of the fractal uncertainty principle \cite[Lemma~3.1]{BD}. Their multiplier construction uses the fractal structure of the set to obtain near-exponential damping at a logarithmic scale excluded by the global one-dimensional logarithmic integral condition. In higher dimensions, Vasilyev's nonradial criterion leaves a gap between the power $d+1-\varepsilon$ for every $\varepsilon>0$ and the insufficient critical power $d+1$. This suggests asking whether the additional power in \eqref{eq:old-power} can be replaced by a logarithmic correction, and what the sharp correction should be.

Vasilyev estimates the maximum on each dyadic annulus. We instead control the exact spherical envelope through the mass forced by a Lipschitz peak on a spherical cap. Combined with Theorem~\ref{thm:main}, this yields our second main result, a sharp logarithmic volume criterion.

For $\beta\in\R$, define
\begin{equation}\label{eq:weighted-integrals}
 I_{\beta}(\Omega)
 =\int_{\R^d}\frac{\Omega(x)[\log(e+|x|)]^\beta}
 {(1+|x|)^{d+1}}\dd x.
\end{equation}
(Here $e$ is chosen only to avoid a singularity at the origin for arbitrary real $\beta$. Since we are concerned with growth at infinity, it may be replaced by any fixed constant $a>1$. Similar choices below are made for the same reason.)
\begin{theorem}\label{thm:log}
Let $d\ge2$ and let $\Omega:\R^d\to[0,\infty)$ be globally Lipschitz.
If
\begin{equation}\label{eq:universal}
 \int_{\R^d}
 \frac{\Omega(x)[\log(e+|x|)]^\beta}{(1+|x|)^{d+1}}\,\dd x
 <\infty
 \qquad\text{for some }\beta>d-1,
\end{equation}
then $e^{-\Omega}$ admits radial BM multipliers.

The threshold $\beta>d-1$ is sharp. For every $\beta\le d-1$, there exists a nonnegative globally Lipschitz
$\Omega\in C^\infty(\R^d)$ with $I_\beta(\Omega)<\infty$ such that every $f\in L^2(\R^d)$ with compact Fourier support and $|f|\le e^{-\Omega}$ is identically zero. The counterexample can have arbitrarily small positive Lipschitz constant.
\end{theorem}
The following remark clarifies that this theorem indeed improves Vasilyev's criterion.
\begin{remark}
For every $\varepsilon>0$ and $\beta>d-1$, $[\log(e+|x|)]^\beta\lesssim_{\beta,\varepsilon}(1+|x|)^\varepsilon$. Thus \eqref{eq:old-power} implies $I_\beta(\Omega)<\infty$. The weakening is strict. Indeed, fix $p>\beta+1$, and let
$h:[0,\infty)\to[0,\infty)$ be globally Lipschitz with $h(r)=r/(\log r)^p$ for all sufficiently large $r$. Define
\[
 \Omega(x)=h(|x|)
 \left(1+\frac{x_1}{2\sqrt{1+|x|^2}}\right).
\]
It is straightforward to verify that $\Omega$ is nonnegative, nonradial, and globally Lipschitz, with
$I_\beta(\Omega)<\infty$, whereas \eqref{eq:old-power} fails for every $\varepsilon>0$. As already noted by Vasilyev in the discussion following \cite[Theorem~2]{Vas}, the uncorrected integral $I_0$ is
insufficient in dimensions $d\ge2$.
\end{remark}
Finally, we return to Cohen's directional integrability condition. Although this condition does not guarantee the existence of radial multipliers with the prescribed decay, we show that such multipliers exist after a logarithmic loss under Lipschitz regularity alone. The threshold for this loss is sharp.
\begin{theorem}
\label{thm:loss}
Let $d\ge2$, and let $\Omega:\R^d\to[0,\infty)$ be globally Lipschitz with $\mathcal C_{\rm dir}(\Omega)<\infty$.
For every $\gamma>1/2$ and every $\sigma>0$, there exists a nonzero radial function $f\in L^2(\R^d)$ such that
\begin{equation}\label{eq:loss}
 \supp\widehat f\subset\overline B(0,\sigma),\qquad
 |f(x)|\le
 \exp\!\left(-\frac{\Omega(x)}
 {[\log(e+|x|)]^\gamma}\right)
 \quad(x\in\R^d).
\end{equation}

The threshold is sharp: There exists a nonnegative globally Lipschitz function $\Omega\in C^\infty(\R^d)$, vanishing near the origin and satisfying $\mathcal C_{\rm dir}(\Omega)<\infty$, such that, for every $0\leq\gamma\le1/2$ and every $c,C>0$, no nonzero radial bandlimited function $f\in L^2(\R^d)$ satisfies
\begin{equation}\label{eq:sharp-upper}
 |f(x)|\le C\exp\!\left(
 -\frac{c\Omega(x)}{[\log(e+|x|)]^\gamma}
 \right)
 \quad(x\in\R^d).
\end{equation}
\end{theorem}
The sharpness assertion concerns radial multipliers and does not exclude the existence of nonradial multipliers with stronger decay. Constructing nonradial multipliers at $\gamma=0$ under the same assumptions would establish the existence part of Cohen's theorem under Lipschitz regularity alone; obtaining his quantitative bounds as well would resolve the open question raised in Remark \ref{rem:main}.

Another contribution of this paper is a more direct transform proof of Vasilyev's radial theorem \cite{Vas}, stated above as Theorem~\ref{thm:vas-radial}. Recently, Bergman also gave a simplified proof by directly lifting an even entire function to several variables, and obtained a radial Cartwright theorem \cite{Bergman}. Our approach uses the Hankel--Paley--Wiener theorem to transfer compact support of the cosine transform to the Fourier--Bessel transform for every real order $\nu\ge-1/2$, thereby avoiding separate arguments for odd and even dimensions.

\medskip
\noindent\textbf{Organization of the paper.}
In Section~\ref{sec:radial-characterization}, we give a
Hankel--Paley--Wiener proof of Vasilyev's radial theorem and
establish the characterization in Theorem~\ref{thm:main}.
Section~\ref{sec:logarithmic-criterion} proves the sharp logarithmic
volume criterion in Theorem~\ref{thm:log}.
We first derive a general weighted sufficient condition from a
spherical maximum estimate, and then construct tube weights to
prove sharpness, including the double-logarithmic endpoint.
In Section~\ref{sec:cohen-comparison}, we compare the spherical
maximum condition with Cohen's directional integrability
condition and establish the optimal logarithmic loss, yielding
Theorem~\ref{thm:loss}.

\medskip
\noindent\textbf{Notation.}
We use $|x|$ for the Euclidean norm,
$B(x,r)$ for the open ball of radius $r$ centered at $x$,
and $\Sph=\{\theta\in\R^d:|\theta|=1\}$ for the unit sphere.
The notation $\mathrm dS$ denotes its usual surface measure.
We write
\[
 \langle x\rangle=(1+|x|^2)^{1/2},
 \qquad t_+=\max\{t,0\}.
\]
The space $\Lip(\R^d)$ consists of globally Lipschitz functions,
and $\Lip(u)$ denotes the least Lipschitz constant of $u$.
For nonnegative quantities $A$ and $B$, the notation
$A\lesssim B$ means $A\le CB$ for a constant independent of
the variables under consideration; dependence on additional
parameters is indicated by subscripts when needed.
We write $A\asymp B$ if both $A\lesssim B$ and $B\lesssim A$.
Constants denoted by $C$ may change from line to line.
\section{The radial characterization}\label{sec:radial-characterization}

In this section, we prove Theorem~\ref{thm:main} using Vasilyev's radial theorem, Theorem~\ref{thm:vas-radial}. Before doing so, we give a simplified proof of the latter. The idea arose during the author's joint work on the fractal uncertainty principle for the Fourier--Bessel transform \cite{ZWD}: the Hankel--Paley--Wiener theorem provides a direct way to transfer spectral support, thereby avoiding the separate arguments for odd and even dimensions in Vasilyev's proof, including the use of Sonin's formula in even dimensions. For another simplified proof, based on lifting an even entire function to several variables, see Bergman~\cite{Bergman}.

The key step in our proof is to establish the support transfer from the cosine transform to the Fourier--Bessel transform. To this end, for $\nu\ge-1/2$, we introduce the measure
$\mathrm d\mu_\nu(r)=r^{2\nu+1}\mathrm dr$ and the Fourier--Bessel transform
\begin{equation}\label{eq:bessel-transform}
 \mathcal B_\nu g(\rho)
 =\rho^{-\nu}\int_0^\infty
 g(r)J_\nu(r\rho)r^{\nu+1}\dd r,
 \qquad \rho>0.
\end{equation}
Initially defined for $g\in C_c^\infty(0,\infty)$, this transform extends to an isometric isomorphism of $L^2(\mu_\nu)$ onto itself, with
\[
 \mathcal B_\nu^{-1}=\mathcal B_\nu.
\]
For these standard facts, see, for example, Andersen~\cite{N}. We also write
\[
 Tg(t)=\int_0^\infty g(r)\cos(rt)\dd r
\]
for the cosine transform, interpreted by its $L^2$ extension when necessary.

\begin{lemma}\label{lem:cosine-bessel}
Let $\nu\ge-1/2$, $A>0$, and $g\in L^2(0,\infty)\cap L^2(\mu_\nu)$. If $\supp(Tg)\subset[0,A]$, then
\begin{equation}\label{eq:cosine-bessel-support}
 \supp(\mathcal B_\nu g)\subset[0,A].
\end{equation}
Both support assertions are understood in the corresponding $L^2$ spaces.
\end{lemma}
\begin{proof}
Write $h=Tg$. Cosine inversion gives
\[
 g(r)=\frac2\pi\int_0^A h(t)\cos(rt)\dd t
\]
in $L^2(0,\infty)$. Since $h\in L^2(0,A)\subset L^1(0,A)$, the function
\[
 G(z)=\frac2\pi\int_0^A h(t)\cos(zt)\dd t
\]
is even and entire, agrees with $g$ almost everywhere on $(0,\infty)$, and satisfies
\[
 |G(x+iy)|\le\frac2\pi\|h\|_{L^1(0,A)}e^{A|y|}.
\]
Thus $G$ has exponential type at most $A$. Since $G=g$ almost everywhere on $(0,\infty)$ and $g\in L^2(\mu_\nu)$ by hypothesis, the restriction of $G$ to the positive real axis belongs to $L^2(\mu_\nu)$. The Hankel--Paley--Wiener theorem of Griffith \cite{Griffith}, in the form restated in \cite[Theorem 1.1]{Khats}
and rescaled from bandwidth $1$ to $A$, provides a function $q\in L^2((0,A),\mu_\nu)$ such that
\[
 G(r)=r^{-\nu}\int_0^A q(s)J_\nu(rs)s^{\nu+1}\dd s.
\]
For comparison with the unweighted formulation of that theorem, its coefficient is $s^{\nu+1/2}q(s)\in L^2(0,A)$. After extending $q$ by zero, the displayed representation is $g=\mathcal B_\nu q$. The Hankel inversion theorem now gives $\mathcal B_\nu g=q$, which proves \eqref{eq:cosine-bessel-support}. At $\nu=-1/2$ this also follows directly from $\mathcal B_{-1/2}=\sqrt{2/\pi}\,T$.
\end{proof}

For $d\ge1$, $\nu=d/2-1$, and a radial $L^2$ function $F(x)=g(|x|)$, the radial Fourier formula in our normalization reads
\begin{equation}\label{eq:radial-fourier-bessel}
 \begin{aligned}
 \widehat F(\xi)
 &=2\pi|\xi|^{1-d/2}
 \int_0^\infty g(r)J_{d/2-1}(2\pi r|\xi|)r^{d/2}\dd r\\
 &=(2\pi)^{d/2}\mathcal B_{d/2-1}g(2\pi|\xi|).
 \end{aligned}
\end{equation}
For general $L^2$ profiles this identity follows by truncation and Plancherel, without assuming pointwise convergence of the integral. Consequently, Lemma~\ref{lem:cosine-bessel} implies
\[
 \supp(Tg)\subset[0,A]
 \quad\Longrightarrow\quad
 \supp\widehat F\subset\overline B(0,A/(2\pi)),
\]
provided $g\in L^2(0,\infty)\cap L^2(r^{d-1}\mathrm dr)$.

We now give a direct proof of Theorem~\ref{thm:vas-radial}. In this proof we denote the radial logarithmic weight by $\Phi$.

\begin{proof}[Direct proof using Lemma \ref {lem:cosine-bessel} and Theorem \ref{thm:classical-bm}]
Write $\Phi(x)=m(|x|)$ and fix $\sigma>0$. Polar coordinates and the assumed integrability condition give
$\int_\R m(|t|)(1+t^2)^{-1}\dd t<\infty$. Choose an integer $N>d$ and apply the one-dimensional BM Theorem \ref{thm:classical-bm} \cite{BM} to the even nonnegative Lipschitz weight
\[
 Q(t)=\tfrac12m(|t|)+N\log(1+t^2).
\]
The added term preserves logarithmic integrability, since $\log(1+t^2)/(1+t^2)\in L^1(\R)$. It yields $0\ne h\in L^2(\R)$ with
\[
 \supp\widehat h\subset[-\sigma/2,\sigma/2],
 \qquad |h(t)|\le e^{-m(|t|)/2}(1+t^2)^{-N}.
\]
Set $q(t)=h(t)h(-t)$. The bandlimited function $h$ extends to an entire function, so $q$ is a nonzero even entire function. Moreover,
\[
 |q(t)|\le e^{-m(|t|)}(1+t^2)^{-2N},
 \qquad \supp\widehat q\subset[-\sigma,\sigma].
\]
The support assertion follows from the product--convolution identity: $\widehat h$ is compactly supported and belongs to $L^1\cap L^2$, and reflection preserves its support interval. Let $g=q|_{[0,\infty)}$. The decay bound gives
\[
 g\in L^2(0,\infty)\cap L^2((0,\infty),r^{d-1}\mathrm dr).
\]
Evenness gives $\widehat q(\xi)=2Tg(2\pi\xi)$ for $\xi>0$, hence $\supp(Tg)\subset[0,2\pi\sigma]$. Apply Lemma~\ref{lem:cosine-bessel} with $\nu=d/2-1$ and \eqref{eq:radial-fourier-bessel}. The function $F(x)=g(|x|)$ is in $L^2(\R^d)$, has Fourier support in $\overline B(0,\sigma)$, and satisfies $|F(x)|\le e^{-\Phi(x)}$.
It is nonzero because the nonzero even entire function $q$ cannot vanish on the whole positive real axis. This proves the radial theorem in every dimension by the same argument.
\end{proof}
\begin{remark}
The classical Hankel--Paley--Wiener theorem yields the support conclusion in every dimension. The proof establishes the qualitative radial BM theorem. The quantitative local lower bounds in \cite[Proposition 1]{Vas} involve a separate one-dimensional input.
\end{remark}

We now proceed to prove Theorem~\ref{thm:main}. For a nonnegative, locally bounded radial function
$\Phi(x)=m(|x|)$, the integral condition in Theorem~\ref{thm:vas-radial} is equivalent to
\begin{equation}\label{eq:radial-integral}
 \int_1^\infty\frac{m(r)}{r^2}\dd r<\infty.
\end{equation}
Indeed, polar coordinates give the weight $r^{d-1}(1+r^2)^{-(d+1)/2}$, which is comparable to $r^{-2}$ for $r\ge1$; the integral near zero is finite by local boundedness. Likewise,
\[
 \mathcal J(\Phi)
 =\int_0^\infty\frac{m(r)}{1+r^2}\dd r<\infty
 \quad\Longleftrightarrow\quad
 \int_1^\infty\frac{m(r)}{r^2}\dd r<\infty.
\]
Throughout the paper, whenever only convergence or divergence is at issue, we use the equivalent tail integral in
\eqref{eq:radial-integral} for convenience, without further comment.

\begin{proof}[Proof of Theorem~\ref{thm:main}]
If $L$ is a Lipschitz constant for $\Omega$, then
\[
 |m_\Omega(r)-m_\Omega(s)|\le L|r-s|,
 \qquad r,s\ge0.
\]
Indeed, $\Omega(r\theta)\le\Omega(s\theta)+L|r-s|$ for each $\theta$; take maxima and interchange $r$ and $s$.
Thus $\Phi(x)=m_\Omega(|x|)$ is a nonnegative globally Lipschitz radial majorant of $\Omega$. If \textup{(i)} holds, local boundedness and \eqref{eq:radial-integral} show that $\Phi$ satisfies Theorem~\ref{thm:vas-radial}. For each $\sigma>0$, that theorem gives a nonzero radial $f_\sigma$ with the required support and $|f_\sigma|\le e^{-\Phi}\le e^{-\Omega}$. This proves \textup{(i)}$\Rightarrow$\textup{(iii)}, and \textup{(iii)}$\Rightarrow$\textup{(ii)} is immediate.

Assume \textup{(ii)}. Fix $x\in\R^d$ and put $r=|x|$. By continuity of $\Omega$ and compactness of the sphere of radius $r$, there is a point $y$ with $|y|=r$ such that $\Omega(y)=m_\Omega(r)$. Since $f$ is radial and
$|x|=|y|$, we have $f(x)=f(y)$. Applying the assumed bound at $y$ gives
\[
 |f(x)|=|f(y)|
 \le e^{-\Omega(y)}
 =e^{-m_\Omega(r)}
 =e^{-\Phi(x)}.
\] 
Applying Bergman's necessary condition \cite[Proposition~1]{Bergman} to the radial logarithmic weight $\Phi$
and the nonzero bandlimited function $f$ gives
\[
 \int_{\R^d}\frac{\Phi(x)}{(1+|x|^2)^{(d+1)/2}}\dd x<\infty.
\]
By \eqref{eq:radial-integral}, this is equivalent to $\mathcal J(\Omega)<\infty$, which proves \textup{(i)}.
\end{proof}
\section{Proof of the sharp logarithmic criterion}
\label{sec:logarithmic-criterion}

In this section, we prove Theorem~\ref{thm:log} in two parts. We first establish sufficiency, using a spherical maximum estimate to derive a multiplier criterion with a general integral weight. We then construct counterexamples to prove sharpness.

\subsection{The spherical maximum estimate}

 We next show that the spherical maximum is controlled by the spherical mass. This estimate is a special case of Gagliardo--Nirenberg interpolation on the sphere; for a formulation on general compact hypersurfaces, see \cite[Theorem~1\textup{(vi)}, equation~(2.17)]{DDM}. We include an elementary spherical cap proof to make the dependence
on $r$ and $L$ explicit.

\begin{lemma}\label{lem:cap}
Let $u:\R^d\to[0,\infty)$ be globally $L$-Lipschitz, where $L>0$, and assume $u(0)=0$. Define, for $r\ge0$,
\[
 M(r)=\max_{\theta\in\Sph}u(r\theta),
 \qquad A(r)=\int_{\Sph}u(r\theta)\dd S(\theta).
\]
Then $M$ is $L$-Lipschitz on $[0,\infty)$, $0\le M(r)\le Lr$, and
\begin{equation}\label{eq:cap-bound}
 M(r)^d\le C_d(Lr)^{d-1}A(r)\qquad(r>0).
\end{equation}
Moreover $x\mapsto M(|x|)$ is a globally $L$-Lipschitz majorant of $u$.
\end{lemma}

\begin{proof}
The Lipschitz continuity of $M$ follows from the argument in the proof of Theorem~\ref{thm:main}. The bound $M(r)\le Lr$ follows from $u(0)=0$.

Fix $r>0$. If $M(r)=0$ there is nothing to prove. Otherwise choose $\theta_0$ attaining the maximum and put
\[
 \delta=\frac{M(r)}{2Lr}\in(0,1/2].
\]
On the cap $\{\theta\in\Sph:|\theta-\theta_0|\le\delta\}$, the Lipschitz condition gives
\[
 u(r\theta)\ge M(r)-Lr|\theta-\theta_0|\ge\frac12M(r).
\]
This cap has surface measure at least $c_d\delta^{d-1}$. For completeness, rotate $\theta_0$ to $e_1$ and parametrize the cap by $\theta=(\sqrt{1-|v|^2},v)$ with $|v|\le\delta/2$ in $\R^{d-1}$. Its chordal distance to $e_1$ is at most $\delta$, and the area Jacobian is $(1-|v|^2)^{-1/2}\ge1$. Thus
\[
 A(r)\ge c_d M(r)\delta^{d-1}
       =c'_d\frac{M(r)^d}{(Lr)^{d-1}},
\]
which proves \eqref{eq:cap-bound}. Finally, $u(x)\le M(|x|)$ and $|M(|x|)-M(|y|)|\le L\,\big||x|-|y|\big|\le L|x-y|$.
\end{proof}
We now use the spherical maximum estimate to derive a multiplier criterion with a general integral weight.
\begin{proposition}\label{prop:general}
Let $d\ge2$, and let $W:[1,\infty)\to(0,\infty)$ be measurable. Assume
\begin{equation}\label{eq:J}
 J_W:=\int_1^\infty W(t)^{-1/(d-1)}\dd t<\infty.
\end{equation}
If $\Omega:\R^d\to[0,\infty)$ is globally Lipschitz and
\begin{equation}\label{eq:K}
 K_W(\Omega):=
 \int_{|x|\ge e}\frac{\Omega(x)W(\log|x|)}{|x|^{d+1}}\dd x<\infty,
\end{equation}
then $\mathcal J(\Omega)<\infty$. Consequently $e^{-\Omega}$ admits radial BM multipliers by Theorem~\ref{thm:main}.
\end{proposition}

\begin{proof}
We first reduce to the case of a weight vanishing at the origin, as required by Lemma~\ref{lem:cap}. Set $c_0=\Omega(0)$ and $u=(\Omega-c_0)_+$. Then $u$ is nonnegative and globally Lipschitz, with $u(0)=0$ and
\[
 0\le u\le\Omega\le c_0+u.
\]
The first comparison gives $K_W(u)\le K_W(\Omega)<\infty$. The second gives $m_\Omega(r)\le c_0+m_u(r)$ for every $r\ge0$, so
\[
 \mathcal J(\Omega)
 \le c_0\int_0^\infty\frac{\dd r}{1+r^2}+\mathcal J(u)
 = \frac{\pi}{2}c_0+\mathcal J(u).
\]
It therefore suffices to prove that $\mathcal J(u)<\infty$. Choose a positive Lipschitz constant $L$ for $u$ and let $M$ and $A$ be its spherical maximum and spherical mass, respectively, as in Lemma~\ref{lem:cap}.

Polar coordinates and \eqref{eq:cap-bound} yield
\begin{align}
 \int_e^\infty\frac{M(r)^d W(\log r)}{r^{d+1}}\,\dd r
 &\le C_d L^{d-1}\int_e^\infty
       \frac{A(r)W(\log r)}{r^2}\,\dd r \notag\\
 &= C_d L^{d-1}\int_e^\infty\int_{\Sph}
       \frac{u(r\theta)W(\log r)}{r^{d+1}}\,
       r^{d-1}\,\dd S(\theta)\,\dd r \notag\\
 &= C_d L^{d-1}\int_{|x|\ge e}
       \frac{u(x)W(\log|x|)}{|x|^{d+1}}\,\dd x \notag\\
 &\le C_d L^{d-1}K_W(\Omega).
 \label{eq:energy-M}
\end{align}
Write $r=e^t$ and $b(t)=e^{-t}M(e^t)$. Then
\[
 \int_e^\infty\frac{M(r)^dW(\log r)}{r^{d+1}}\dd r=\int_1^\infty b(t)^dW(t)\dd t.
\]
H\"older's inequality with exponents $d$ and $d/(d-1)$ now gives
\begin{align}
 \int_e^\infty\frac{M(r)}{r^2}\dd r
 &=\int_1^\infty b(t)\dd t\notag\\
 &\le
 \left(\int_1^\infty b(t)^dW(t)\dd t\right)^{1/d}
 J_W^{(d-1)/d}\notag\\
 &\le C_d L^{(d-1)/d}K_W(\Omega)^{1/d}
 J_W^{(d-1)/d}<\infty.\label{eq:envelope-estimate}
\end{align}
Thus $\mathcal J(u)<\infty$, Theorem~\ref{thm:main} gives the asserted radial multipliers.
\end{proof}

\begin{proof}[Sufficiency in Theorem~\ref{thm:log}]
Take $W(t)=t^\beta$. For $|x|\ge e$, the integrand in \eqref{eq:K} is comparable to that in $I_\beta$; hence $K_W<\infty$. Also
\[
 J_W=\int_1^\infty t^{-\beta/(d-1)}\dd t<\infty
 \quad\Longleftrightarrow\quad\beta>d-1.
\]
Proposition~\ref{prop:general} proves \eqref{eq:universal}.
\end{proof}

\subsection{A tube construction and sharpness}
In this subsection, we construct counterexamples to prove sharpness. As in the sufficiency argument, we first establish a general principle and then apply it to specific examples. To this end, we begin by recalling a classical logarithmic integral theorem.

The following lemma is a direct consequence of the classical logarithmic integral theorem for bounded holomorphic functions in the upper half-plane \cite[Chapter~III, pp.~49--51]{Koosis}. See also \cite{HavinJoricke} for the role of logarithmic integrals in uncertainty principles. For completeness, we recall a standard proof based on Jensen's formula \cite[Chapter~I]{Koosis}. Throughout, following \cite{Koosis}, we write $\log^-s=\max\{0,-\log s\}$ for $s>0$, with $\log^-0=+\infty$.

\begin{lemma}\label{lem:one-dimensional-necessity}
Let $g$ be a nonzero entire function satisfying
\[
 |g(z)|\le A e^{S|\operatorname{Im}z|}\quad(z\in\C)
\]
for some $A>0$ and $S\ge0$. Then
\begin{equation}\label{eq:log-necessity}
 \int_\R\frac{\log^-|g(t)|}{1+t^2}\dd t<\infty.
\end{equation}
Consequently, if $q\ge0$, $|g(t)|\le e^{-q(t)}$, and $\int_\R q(t)(1+t^2)^{-1}\dd t=\infty$, then $g\equiv0$.
\end{lemma}

\begin{proof}
To remove the exponential growth in the upper half-plane, set $H(z)=e^{iSz}g(z)$. For $z=x+iy$ with $y>0$, we have
\[
 |H(z)|=e^{-Sy}|g(z)|\le A.
\]
Moreover, $H\not\equiv0$ and $|H(t)|=|g(t)|$ for $t\in\R$.

Choose $z_0=a+ib$, with $b>0$, such that $H(z_0)\ne0$. To apply Jensen's formula at this nonzero value, use the conformal map
\[
 \varphi(w)=a+ib\frac{1+w}{1-w},
 \qquad |w|<1,
\]
which maps the unit disk onto the upper half-plane and satisfies $\varphi(0)=z_0$. Indeed,
\[
 \operatorname{Im}\varphi(w)
 =b\frac{1-|w|^2}{|1-w|^2}>0.
\]
Put $A_1=\max(A,1)$ and
\[
 F(w)=\frac{H(\varphi(w))}{A_1}.
\]
Then $F$ is holomorphic in the unit disk, $|F|\le1$, and $F(0)\ne0$. In particular, $-\log|F|$ is nonnegative. Jensen's formula, or the submean inequality for $\log|F|$ (see \cite[Chapter~I]{Koosis}), gives
\[
 \frac1{2\pi}\int_0^{2\pi}
 -\log|F(\rho e^{i\theta})|\,\dd\theta
 \le-\log|F(0)|<\infty,
 \qquad 0<\rho<1.
\]
Thus these nonnegative integrals have a finite upper bound independent of $\rho$.

For $0<\theta<2\pi$, the boundary value of the conformal map is
\[
 \varphi(e^{i\theta})
 =a-b\cot\frac{\theta}{2}\in\R.
\]
Continuity of $H$ up to the real axis therefore gives
\[
 \lim_{\rho\uparrow1}F(\rho e^{i\theta})
 =\frac{H(\varphi(e^{i\theta}))}{A_1}.
\]
At a zero of this boundary value, the corresponding negative logarithm tends to $+\infty$. The only exceptional boundary point is $w=1$, which is mapped to infinity and has angular measure zero. Fatou's lemma applied to the nonnegative functions $-\log|F(\rho e^{i\theta})|$ now yields
\[
 \frac1{2\pi}\int_0^{2\pi}
 -\log\frac{|H(\varphi(e^{i\theta}))|}{A_1}\,\dd\theta
 \le-\log|F(0)|<\infty.
\]
Changing variables back to the real line and using $|H(t)|=|g(t)|$, $A_1\ge1$, and the comparability of the resulting Poisson weight with $(1+t^2)^{-1}$, we obtain
\[
 \int_\R\frac{\log^-|g(t)|}{1+t^2}\,\dd t<\infty,
\]
which proves \eqref{eq:log-necessity}.

Finally, if $q\ge0$ and $|g(t)|\le e^{-q(t)}$, then $\log^-|g(t)|\ge q(t)$. Divergence of
$\int_\R q(t)(1+t^2)^{-1}\,\dd t$ therefore contradicts \eqref{eq:log-necessity} whenever $g\not\equiv0$.
\end{proof}
We next apply Lemma~\ref{lem:one-dimensional-necessity} to show that a majorant admits no BM multiplier if its logarithmic integral diverges along every line parallel to a fixed direction.
\begin{lemma}\label{lem:parallel-lines}
Let $\Omega:\R\times\R^{d-1}\to[0,\infty)$ be continuous and satisfy
\begin{equation}\label{eq:all-lines}
 \int_\R\frac{\Omega(t,y)}{1+t^2}\dd t=\infty
 \qquad\text{for every fixed }y\in\R^{d-1}.
\end{equation}
Then no nonzero compact-spectrum function $f\in L^2(\R^d)$ can satisfy $|f|\le e^{-\Omega}$ almost everywhere.
\end{lemma}

\begin{proof}
If $\supp\widehat f\subset\overline B(0,S)$, then $\widehat f\in L^1$ by Cauchy--Schwarz on a finite-measure ball. Hence
\[
 F(z)=\int_{\R^d}\widehat f(\xi)e^{2\pi iz\cdot\xi}\dd\xi,
 \qquad z\in\C^d,
\]
is entire, represents $f$ continuously on $\R^d$, and obeys $|F(x+i\eta)|\le C_f e^{2\pi S|\eta|}$. The almost-everywhere upper bound extends to every real point by continuity.

For each fixed $y\in\R^{d-1}$, the function $g_y(z)=F(z,y)$ satisfies the growth hypothesis of Lemma~\ref{lem:one-dimensional-necessity} and $|g_y(t)|\le e^{-\Omega(t,y)}$. By \eqref{eq:all-lines}, it is identically zero. This holds for every $y$, so $f\equiv0$.
\end{proof}

We next construct a nonnegative weight supported in a tube containing the tails of every line parallel to the $x_1$-axis. A suitable choice of the profile will make its Poisson integral diverge along each of these lines, so that Lemma~\ref{lem:parallel-lines} excludes every nonzero bandlimited multiplier.
Write $x=(x_1,x')\in\R\times\R^{d-1}$. The following lemma also quantifies the angular concentration of this
construction.

\begin{lemma}\label{lem:tube}
Let $g:[0,\infty)\to[0,\infty)$ be globally Lipschitz and, for all sufficiently large $r$, assume
\[
 0<g(r)\le r/4,\qquad g(r)\longrightarrow\infty,
 \qquad g\text{ is nondecreasing}.
\]
Define
\begin{equation}\label{eq:tube}
 \Omega_g(x)=\big(g(|x|)-|x'|\big)_+.
\end{equation}
Then $\Omega_g$ is nonnegative and globally Lipschitz. For large $r$,
\begin{equation}\label{eq:tube-average}
 c_d\frac{g(r)^d}{r^{d-1}}
 \le\int_{\Sph}\Omega_g(r\theta)\dd S(\theta)
 \le C_d\frac{g(r)^d}{r^{d-1}}.
\end{equation}
If $\int_1^\infty g(r)r^{-2}\dd r=\infty$, then $\Omega_g$ satisfies
\eqref{eq:all-lines}.
\end{lemma}

\begin{proof}
The  Lipschitz of $\Omega_g$ follows from the fact that: the maps $x\mapsto g(|x|)$ and $x\mapsto|x'|$ have Lipschitz constants at most $\Lip(g)$ and $1$, respectively, and taking positive part does not increase the constant. 

To prove \eqref{eq:tube-average}, fix $r$ sufficiently large that $0<g(r)\le r/4$, and put
$\delta=g(r)/r\in(0,1/4]$. Writing $\theta=(\theta_1,\theta')\in\Sph$, we have
\[
 \Omega_g(r\theta)=\bigl(g(r)-r|\theta'|\bigr)_+.
\]
Thus $\Omega_g(r\theta)>0$ precisely when $|\theta'|<\delta$. This region consists of two caps around $e_1$ and $-e_1$, parametrized by
\[
 \theta_\pm(y)=\bigl(\pm\sqrt{1-|y|^2},y\bigr),
 \qquad y\in\R^{d-1},\quad |y|<\delta.
\]
For either parametrization, the surface element is
\[
 \dd S(\theta_\pm(y))
 =\frac{\dd y}{\sqrt{1-|y|^2}}.
\]
Since $\delta\le1/4$, we have
\[
 1\le\frac1{\sqrt{1-|y|^2}}\le2
 \qquad (|y|\le\delta).
\]
Consequently, the exact integral representation
\[
 \int_{\Sph}\Omega_g(r\theta)\,\dd S(\theta)
 =2\int_{|y|<\delta}
   \frac{g(r)-r|y|}{\sqrt{1-|y|^2}}\,\dd y
\]
gives the upper bound
\[
 \int_{\Sph}\Omega_g(r\theta)\,\dd S(\theta)
 \le4g(r)\int_{|y|<\delta}\dd y
 \le C_d g(r)\delta^{d-1}
 =C_d\frac{g(r)^d}{r^{d-1}}.
\]
For the lower bound, restrict the integral to $|y|\le\delta/2$. On this smaller ball,
\[
 g(r)-r|y|
 \ge g(r)-\frac{r\delta}{2}
 =\frac{g(r)}2.
\]
Using also $(1-|y|^2)^{-1/2}\ge1$, we obtain
\[
 \int_{\Sph}\Omega_g(r\theta)\,\dd S(\theta)
 \ge 2\cdot\frac{g(r)}2
       \int_{|y|\le\delta/2}\dd y
 \ge c_d g(r)\delta^{d-1}
 =c_d\frac{g(r)^d}{r^{d-1}}.
\]
This proves \eqref{eq:tube-average}.

Fix $x'\in\R^{d-1}$. Since $g(r)\to\infty$ and $g$ is eventually nondecreasing, there exists $T>0$ such that, for every $x_1\ge T$,
\[
 g(x_1)\ge2|x'|,
 \qquad
 g\bigl(\sqrt{x_1^2+|x'|^2}\bigr)\ge g(x_1).
\]
Consequently,
\[
 \Omega_g(x_1,x')
 =\bigl(g(\sqrt{x_1^2+|x'|^2})-|x'|\bigr)_+
 \ge g(x_1)-|x'|
 \ge\tfrac12g(x_1).
\]
The assumed divergence of $\int_1^\infty g(r)r^{-2}\,\dd r$ therefore gives
\[
 \int_\R\frac{\Omega_g(x_1,x')}{1+x_1^2}\,\dd x_1
 \ge\frac12\int_T^\infty\frac{g(x_1)}{1+x_1^2}\,\dd x_1
 =\infty.
\]
Since $x'$ was arbitrary, this proves \eqref{eq:all-lines}. The divergent positive tail already implies \eqref{eq:all-lines}.
\end{proof}
Now, we use Lemma \ref{lem:tube} to prove the necessity part in Theorem~\ref{thm:log}.
\begin{proof}[Necessity and endpoint failure in Theorem~\ref{thm:log}]
Choose a number $a$ with
\begin{equation}\label{eq:a-range}
 \frac1d<a\le1.
\end{equation}
For large $r$, set
\begin{equation}\label{eq:single-profile}
 g(r)=\frac{r}{\log r\,(\log\log r)^a}.
\end{equation}
Since $g'(r)\to0$ as $r\to\infty$, we have $|g'(r)|\le1$ for all sufficiently large $r$, so $g$ is $1$-Lipschitz on $[R,\infty)$ for sufficiently large $R$. Extend $g$ nonnegatively and globally Lipschitz on the remaining compact interval. For example, choose a smooth cutoff that is zero below a large $R$ and one above $2R$, multiply the displayed profile by that cutoff where the logarithms are defined, and set it to zero below $R$.
For large $r$,
\[
 g'(r)=\frac{1}{\log r(\log\log r)^a}
 \left(1-\frac1{\log r}
 -\frac{a}{\log r\log\log r}\right)>0.
\]
Thus $g$ meets the hypotheses of Lemma~\ref{lem:tube}, and
\begin{equation}\label{eq:line-diverges}
 \int_{R}^\infty\frac{g(r)}{r^2}\dd r
 =\int_{R}^\infty\frac{\dd r}{r\log r(\log\log r)^a}
 =\infty
\end{equation}
because $a\le1$.

On the other hand, choose $R$ sufficiently large that \eqref{eq:tube-average} and \eqref{eq:single-profile} hold for $r\ge R$. Starting from the definition of $I_\beta$ and using polar coordinates, we obtain
\begin{align}
 I_\beta(\Omega_g)
 &=\int_{\R^d}
   \frac{\Omega_g(x)[\log(e+|x|)]^\beta}
        {(1+|x|)^{d+1}}\,\dd x \notag\\
 &=C_0+\int_R^\infty
   \frac{[\log(e+r)]^\beta r^{d-1}}{(1+r)^{d+1}}
   \left(\int_{\Sph}\Omega_g(r\theta)\,\dd S(\theta)\right)\dd r
   \notag\\
 &\le C_0+C_{d,\beta}\int_R^\infty
   \frac{(\log r)^\beta}{r^2}
   \frac{g(r)^d}{r^{d-1}}\,\dd r \notag\\
 &=C_0+C_{d,\beta}\int_R^\infty
   \frac{g(r)^d(\log r)^\beta}{r^{d+1}}\,\dd r \notag\\
 &=C_0+C_{d,\beta}\int_R^\infty
   \frac{\dd r}
   {r(\log r)^{d-\beta}(\log\log r)^{ad}},
 \label{eq:counter-integral}
\end{align}
where $C_0<\infty$ is the contribution from $|x|<R$. In the inequality we used \eqref{eq:tube-average},
$1+r\asymp r$, and $\log(e+r)\asymp\log r$ for $r\ge R$. The last equality follows by substituting \eqref{eq:single-profile}. This integral is finite if $\beta<d-1$. At $\beta=d-1$ it is still finite because $ad>1$. Hence the same $\Omega_g$ works simultaneously for all $\beta\le d-1$. Lemmas~\ref{lem:tube} and
\ref{lem:parallel-lines} show that there is no nonzero compact-spectrum multiplier, for any spectral radius.

To obtain smooth examples, let $\eta\ge0$ be a $C^\infty$ mollifier of integral one, supported in $B(0,1)$, and set $\widetilde\Omega=\eta*\Omega_g$. It is straightforward to verify that $\widetilde\Omega$ is smooth, nonnegative, nonradial, and globally Lipschitz, with $I_\beta(\widetilde\Omega)<\infty$ and \eqref{eq:all-lines}. Multiplying $\widetilde\Omega$ by a sufficiently small positive constant makes its Lipschitz constant arbitrarily small while preserving these properties.
\end{proof}

\begin{remark}
In this remark, we explain why divergence along a single line does not suffice to obtain a counterexample. A multivariable entire function can vanish on one real line without vanishing everywhere. The construction above forces divergence on \emph{every} line parallel to $e_1$, because the absolute tube width $g(r)$ tends to infinity even though its angular width $g(r)/r$ tends to zero. This is the point that makes the counterexample conclusive.
\end{remark}
In fact, Proposition~\ref{prop:general} and Lemma~\ref{lem:tube} allow us to incorporate any finite number of iterated logarithmic factors; in particular, the case of two such factors gives the following sharp double-logarithmic criterion. For $\beta,\gamma\in\R$, define
\[
 I_{\beta,\gamma}(\Omega)
 :=\int_{\R^d}
 \frac{\Omega(x)[\log(e+|x|)]^\beta
 [\log\log(e^e+|x|)]^\gamma}
 {(1+|x|)^{d+1}}\,\dd x.
\]

\begin{corollary}\label{cor:double}
Let $d\ge2$ and let $\Omega:\R^d\to[0,\infty)$ be globally Lipschitz. If $I_{\beta,\gamma}(\Omega)<\infty$ for some $\beta,\gamma\in\R$ satisfying
\begin{equation}\label{eq:double-region}
 \beta>d-1,
 \qquad\text{or}\qquad
 \beta=d-1\ \text{ and }\ \gamma>d-1,
\end{equation}
then $e^{-\Omega}$ admits radial BM multipliers.

The range \eqref{eq:double-region} is sharp. For every pair $(\beta,\gamma)$ outside this range, there exists a nonnegative globally Lipschitz $\Omega\in C^\infty(\R^d)$ with $I_{\beta,\gamma}(\Omega)<\infty$ such that every
$f\in L^2(\R^d)$ with compact Fourier support and $|f|\le e^{-\Omega}$ is identically zero.
\end{corollary}

\begin{proof}
For sufficiency, apply Proposition~\ref{prop:general} with $W(t)=t^\beta[\log(e+t)]^\gamma$. Then $J_W<\infty$ precisely when \eqref{eq:double-region} holds, and the weighted integral $K_W(\Omega)$ is finite whenever $I_{\beta,\gamma}(\Omega)<\infty$.

For necessity, choose $a$ as in \eqref{eq:a-range} and set, for sufficiently large $r$,
\[
 g(r)=\frac{r}
 {\log r\,\log\log r\,(\log\log\log r)^a}.
\]
Extend $g$ nonnegatively and globally Lipschitz to $[0,\infty)$ using a smooth cutoff, and define $\Omega_g$ by \eqref{eq:tube}. The remaining estimates and the smoothing argument follow the same steps as in the proof of Theorem~\ref{thm:log}, yielding the required counterexamples outside the region \eqref{eq:double-region}.
\end{proof}
\section{Comparison with Cohen's directional condition}\label{sec:cohen-comparison}
In this section, we give a quantitative comparison between our integrability condition and Cohen's directional integrability condition, from which the proof of Theorem~\ref{thm:loss} follows directly.

Cohen's Theorem~1.4 \cite{Cohen}, with his logarithmic weight written as $-\Omega$, assumes $\mathcal C_{\rm dir}(\Omega)<\infty$ together with
\[
 |D^\alpha\Omega(x)|\le C_{\rm reg}\langle x\rangle^{1-|\alpha|},
 \qquad 0\le|\alpha|\le3,
\]
and a normalization $\Omega=0$ near the origin. Its conclusion includes a quantitative local lower bound and an upper bound of the form $|f(x)|\le C e^{-c\sigma\Omega(x)}$.

\begin{proposition}\label{prop:cohen-comparison}
If $\Omega\ge0$ is globally Lipschitz and $\mathcal J(\Omega)<\infty$, then $\mathcal C_{\rm dir}(\Omega)<\infty$. More precisely,
\[
 \mathcal C_{\rm dir}(\Omega)
 \le3\int_0^\infty\frac{m_\Omega(t)}{1+t^2}\dd t.
\]
The converse can fail even for smooth weights satisfying Cohen's symbol estimates through order three and vanishing near the origin.
\end{proposition}
\begin{proof}
We divide the proof into two parts.

\medskip
\noindent\emph{Part 1: The implication.}
For every $r>0$,
\[
 \sup_{|x|=r}G_\Omega(x)
 \le\int_{1/2}^{2}m_\Omega(sr)\dd s.
\]
Tonelli's theorem and the change of variables $t=sr$ give
\[
 \mathcal C_{\rm dir}(\Omega)
 \le\int_{1/2}^{2}\int_0^\infty
       \frac{s\,m_\Omega(t)}{s^2+t^2}\dd t\dd s
 \le3\int_0^\infty\frac{m_\Omega(t)}{1+t^2}\dd t,
\]
where we used
\[
 \frac{s}{s^2+t^2}\le\frac2{1+t^2},
 \qquad \frac12\le s\le2.
\]
The last integral is finite by the assumption $\mathcal J(\Omega)<\infty$.

\medskip
\noindent\emph{Part 2: Failure of the converse.}
We next construct an explicit example showing that the converse fails even under Cohen's symbol estimates through order three. Set
\[
 R_j=16^j,\qquad a_j=j^{-1/3},
\]
and choose a nonzero nonnegative $\phi\in C_c^\infty((1,3/2))$ and $\eta\in C_c^\infty((-1,1))$ with $0\le\eta\le1$ and $\eta(0)=1$. For $x\ne0$, put
\begin{equation}\label{eq:moving-peak}
 \Omega_j(x)=R_j a_j^3\,
 \phi\left(\frac{|x|}{R_j}\right)
 \eta\left(
 \frac{x_1/|x|-|x|/R_j+1}{a_j}\right),
 \qquad
 \Omega(x)=\sum_{j\ge1}\Omega_j(x).
\end{equation}
Each summand extends by zero near the origin. Their supports are contained in the pairwise disjoint annuli
$R_j<|x|<3R_j/2$, so the sum is locally finite, smooth, nonnegative, and vanishes for $|x|\le16$.

Differentiation on the $j$th annulus gives, for $0\le k\le3$,
\[
 |D^k\Omega_j(x)|
 \le C_kR_j^{1-k}a_j^{3-k}
 \le C_kR_j^{1-k}.
\]
Indeed, derivatives of order $\ell\ge1$ of the argument of $\eta$ are bounded by $C_\ell a_j^{-1}R_j^{-\ell}$,
while derivatives of $\phi(|x|/R_j)$ satisfy the usual annular bounds. Since $0<a_j\le1$, the product and chain rules yield the displayed estimate. Thus $\Omega$ is globally Lipschitz and satisfies Cohen's symbol estimates through order three.

For $R_j<r<3R_j/2$, a direction with $\theta_1=\frac{x_1}{|x|}=r/R_j-1\in(0,1/2)$ exists because $d\ge2$. Since $\eta(0)=1$ and $0\le\eta\le1$, we obtain
\[
 m_\Omega(r)=R_j a_j^3\phi(r/R_j).
\]
Consequently,
\[
 \mathcal J(\Omega)
 \asymp\left(\int_1^{3/2}\frac{\phi(t)}{t^2}\dd t\right)
   \sum_{j\ge1}a_j^3
 =\infty.
\]

On a fixed ray $sr\theta$, the factor involving $\eta$ is nonzero only if
\[
 \left|\theta_1-\frac{sr}{R_j}+1\right|<a_j.
\]
The admissible $s$ interval therefore has length at most $2a_jR_j/r$. A contribution also requires
$R_j/2<r<3R_j$. Since $\Omega_j\le R_ja_j^3\|\phi\|_\infty$, it follows that
\[
 \sup_{\theta\in\Sph}\int_{1/2}^{2}\Omega_j(sr\theta)\dd s
 \le C R_ja_j^4\,\mathbf1_{[R_j/2,\,3R_j]}(r).
\]
Using the supremum of a sum bounded by the sum of the suprema and applying Tonelli's theorem, we conclude that
\[
 \mathcal C_{\rm dir}(\Omega)
 \le C\sum_{j\ge1}R_ja_j^4
       \int_{R_j/2}^{3R_j}\frac{\dd r}{1+r^2}
 \le C'\sum_{j\ge1}j^{-4/3}<\infty.
\]
This proves the failure of the converse. By Theorem~\ref{thm:main}, $e^{-\Omega}$ has no nonzero
radial bandlimited multiplier.
\end{proof}
\begin{remark}
The example constructed above shows that when the directions of the spherical peaks vary with the radius, each fixed ray passes through the peaks only over short radial intervals. Consequently, $G_\Omega^*(r)$ can be
substantially smaller than $m_\Omega(r)$, allowing Cohen's directional integrability condition to accommodate a strictly larger class of majorants.
\end{remark}
Although $\mathcal C_{\rm dir}(\Omega)<\infty$ does not imply $\mathcal J(\Omega)<\infty$, we show that it does imply a logarithmically corrected version of this condition, with a sharp threshold for the logarithmic exponent.
\begin{proposition}\label{thm:cohen-log-loss}
Let $d\ge2$, and let $\Omega:\R^d\to[0,\infty)$ be globally Lipschitz. If $\mathcal C_{\rm dir}(\Omega)<\infty$, then
\begin{equation}\label{eq:cohen-log-loss}
 \mathcal J\!\left(
 \frac{\Omega(\cdot)}{[\log(e+|\cdot|)]^\gamma}
 \right)<\infty
 \qquad\text{for every }\gamma>\tfrac12.
\end{equation}
The threshold is sharp: there exists a nonnegative globally Lipschitz function $\Omega\in C^\infty(\R^d)$ such that $\mathcal C_{\rm dir}(\Omega)<\infty$, but
\begin{equation}\label{eq:cohen-log-loss-failure}
 \mathcal J\!\left(
 \frac{\Omega(\cdot)}{[\log(e+|\cdot|)]^\gamma}
 \right)=\infty
 \qquad\text{for every }\gamma\le\tfrac12.
\end{equation}
\end{proposition}

\begin{proof}
We first prove sufficiency. Put $c_0=\Omega(0)$ and
\[
 u(x)=(\Omega(x)-c_0)_+.
\]
Then $u$ is nonnegative and globally Lipschitz, $u(0)=0$, and
\[
 u\le\Omega\le c_0+u,
 \qquad
 \mathcal C_{\rm dir}(u)
 \le\mathcal C_{\rm dir}(\Omega)<\infty.
\]
Choose $L>0$ such that $\operatorname{Lip}(u)\le L$, and write
\[
 M(r)=\max_{\theta\in\Sph}u(r\theta).
\]
We claim that
\begin{equation}\label{eq:directional-quadratic-bound}
 M(r)^2\le2Lr\,G_u^*(r),\qquad r>0.
\end{equation}
If $M(r)=0$, this is immediate. Otherwise, choose $\theta_r\in\Sph$ such that $u(r\theta_r)=M(r)$, and put
\[
 \delta=\frac{M(r)}{2Lr}.
\]
Since $u(0)=0$, we have $M(r)\le Lr$, so $0<\delta\le1/2$. For $|s-1|\le\delta$, the Lipschitz condition gives
\[
 u(sr\theta_r)
 \ge M(r)-Lr|s-1|
 \ge\frac12M(r).
\]
As $[1-\delta,1+\delta]\subset[1/2,2]$, it follows that
\[
 G_u^*(r)
 \ge\int_{1-\delta}^{1+\delta}u(sr\theta_r)\dd s
 \ge\delta M(r)
 =\frac{M(r)^2}{2Lr}.
\]
This proves \eqref{eq:directional-quadratic-bound}. Consequently,
\begin{equation}\label{eq:directional-quadratic-integral}
 \int_1^\infty\frac{M(r)^2}{r^3}\dd r
 \le2L\int_1^\infty\frac{G_u^*(r)}{r^2}\dd r
 \le4L\mathcal C_{\rm dir}(u).
\end{equation}

Fix $\gamma>1/2$. By Cauchy--Schwarz,
\begin{align*}
 \int_e^\infty
 \frac{M(r)}{r^2(\log r)^\gamma}\dd r
 &\le
 \left(\int_e^\infty\frac{M(r)^2}{r^3}\dd r\right)^{1/2}
 \left(\int_e^\infty
 \frac{\dd r}{r(\log r)^{2\gamma}}\right)^{1/2}\\
 &\le
 \frac{2\sqrt{L\mathcal C_{\rm dir}(u)}}{\sqrt{2\gamma-1}}
 <\infty.
\end{align*}
Since $m_\Omega(r)\le c_0+M(r)$, we obtain
\[
 \int_e^\infty
 \frac{m_\Omega(r)}{r^2[\log(e+r)]^\gamma}\dd r<\infty.
\]
The integral over $[0,e]$ is finite by continuity. Finally, division by a radial factor commutes with taking the spherical maximum, so
\[
 \mathcal J\!\left(
 \frac{\Omega(\cdot)}{[\log(e+|\cdot|)]^\gamma}
 \right)
 =
 \int_0^\infty
 \frac{m_\Omega(r)}{(1+r^2)[\log(e+r)]^\gamma}\dd r<\infty.
\]
This proves \eqref{eq:cohen-log-loss}.

We next prove sharpness by slightly modifying the construction in the proof of Proposition~\ref{prop:cohen-comparison}. Choose nonnegative functions $\chi\in C_c^\infty((3/2,2))$ and
$\varphi\in C_c^\infty((-1,1))$ such that $0\le\chi,\varphi\le1$, $\chi\not\equiv0$, and $\varphi(0)=1$.
Set
\[
 R_n=16^n,\qquad
 a_n=\frac{\varepsilon}{\sqrt n\,\log(n+1)},\qquad n\ge2,
\]
where $\varepsilon>0$ is chosen so that $a_n\le1/10$, and define
\begin{equation}\label{eq:directional-sharp-example}
 \Omega(x)=\sum_{n\ge2}
 a_nR_n\chi\!\left(\frac{|x|}{R_n}\right)
 \varphi\!\left(\frac{x_1-R_n}{a_nR_n}\right)
 =:\sum_{n\ge2}\Omega_n(x).
\end{equation}
The summands have disjoint supports in annuli escaping to infinity, and $\|\nabla\Omega_n\|_\infty\le C(a_n+1)\le C$. Thus $\Omega$ is smooth, nonnegative, globally Lipschitz, and vanishes near the origin.

Since $d\ge2$, every sphere with $3R_n/2<r<2R_n$ intersects $\{x_1=R_n\}$. Hence
\[
 m_\Omega(r)=a_nR_n\chi(r/R_n)
\]
on that annulus and is zero outside the union of these annuli. Using $\log(e+R_nt)\asymp n$ for $3/2\le t\le2$, we obtain
\begin{align*}
 \mathcal J\!\left(
 \frac{\Omega(\cdot)}{[\log(e+|\cdot|)]^\gamma}
 \right)
 &=
 \sum_{n\ge2}\int_{3R_n/2}^{2R_n}
 \frac{a_nR_n\chi(r/R_n)}
 {(1+r^2)[\log(e+r)]^\gamma}\dd r\\
 &\asymp\sum_{n\ge2}\frac{a_n}{n^\gamma}
 =\varepsilon\sum_{n\ge2}
 \frac{1}{n^{\gamma+1/2}\log(n+1)}.
\end{align*}
This diverges for every $\gamma\le1/2$.

On the other hand, $\Omega_n(sr\theta)\ne0$ requires
\[
 \frac32R_n<sr<2R_n,\qquad
 |sr\theta_1-R_n|<a_nR_n.
\]
These conditions imply $3R_n/4<r<4R_n$ and $\theta_1>(1-a_n)/2\ge9/20$. Thus the admissible $s$ interval has length at most $2a_nR_n/(r\theta_1)\le Ca_n$, giving
\[
 G_{\Omega_n}^*(r)
 \le Ca_n^2R_n\,\mathbf1_{[3R_n/4,\,4R_n]}(r).
\]
Consequently,
\[
 \mathcal C_{\rm dir}(\Omega)
 \le C\sum_{n\ge2}a_n^2R_n
 \int_{3R_n/4}^{4R_n}\frac{\dd r}{1+r^2}
 \le C'\varepsilon^2\sum_{n\ge2}
 \frac{1}{n[\log(n+1)]^2}<\infty.
\]
This proves sharpness.
\end{proof}
The preceding proposition immediately yields the proof of Theorem~\ref{thm:loss}.
\begin{proof}[Proof of Theorem~\ref{thm:loss}]
For $\gamma>1/2$, the function $\Omega(x)/[\log(e+|x|)]^\gamma$ remains globally Lipschitz, since $\Omega$ is globally Lipschitz and $\Omega(x)\le\Omega(0)+\operatorname{Lip}(\Omega)|x|$. The existence assertion therefore follows by combining Proposition~\ref{thm:cohen-log-loss} with Theorem~\ref{thm:main}, and sharpness follows from the counterexample in that proposition and the necessity of the spherical maximum condition for radial multipliers.
\end{proof}
\section*{Disclosure of AI assistance:}
During the preparation of this manuscript, the author used OpenAI's ChatGPT (Sol 5.6 and Astra 6) mainly to assist with exploring candidate counterexamples, language editing, and LaTeX formatting. The author takes full responsibility for the mathematical correctness, originality, and final content of the manuscript.

\noindent \textbf{Authors' addresses:}

\medskip
\noindent Longben Wei, School of Mathematical Sciences, Guizhou Normal University, Guiyang 550025, P.R. China. \\
Email address: \texttt{longbenwei51@gmail.com}


\begin{thebibliography}{100}
\bibitem{N}
Andersen, N.B.: Real Paley--Wiener theorems for the Hankel transform.
\emph{J. Fourier Anal. Appl.} \textbf{12}, 17--25 (2006). \url{https://doi.org/10.1007/s00041-005-4056-3}.

\bibitem{Bergman}
Bergman, A.: A remark on the Beurling--Malliavin theorem in several variables. Preprint, arXiv:2512.07271 (2025).
\url{https://arxiv.org/abs/2512.07271}.

\bibitem{BM}
Beurling, A. and Malliavin, P.: On Fourier transforms of measures with compact support. \emph{Acta Math.} \textbf{107}, 291--309 (1962).

\bibitem{BD}
Bourgain, J. and Dyatlov, S.: Spectral gaps without the pressure condition. \emph{Ann. of Math.} (2) \textbf{187}, 825--867 (2018).

\bibitem{Cohen}
Cohen, A.: Fractal uncertainty in higher dimensions. \emph{Ann. of Math.} (2) \textbf{202}, 265--307 (2025). \url{https://doi.org/10.4007/annals.2025.202.1.4}.

\bibitem{DDM}
Della Corte, S., Diana, A. and Mantegazza, C.: Uniform Sobolev, interpolation and geometric Calder\'on--Zygmund inequalities for graph hypersurfaces. \emph{Note Mat.} \textbf{44}(1), 53--83 (2024). \url{https://doi.org/10.1285/i15900932v44n1p53}.

\bibitem{Griffith}
Griffith, J.L.: Hankel transforms of functions zero outside a finite interval. \emph{J. Proc. Roy. Soc. New South Wales} \textbf{89},
109--115 (1956). \url{https://doi.org/10.5962/p.360686}.

\bibitem{HanSchlag}
Han, R. and Schlag, W.: A higher-dimensional Bourgain--Dyatlov fractal uncertainty principle. \emph{Anal. PDE} \textbf{13}, 813--863 (2020).

\bibitem{HavinJoricke}
Havin, V. and J\"oricke, B.: \emph{The Uncertainty Principle: A Harmonic Analysis Approach}. Ergebnisse der Mathematik und ihrer Grenzgebiete (3), vol.~28. Springer-Verlag, Berlin (1994).

\bibitem{Hormander}
H\"ormander, L.: Convolution equations in convex domains. \emph{Invent. Math.} \textbf{4}, 306--317 (1968).

\bibitem{JayeMitkovski}
Jaye, B. and Mitkovski, M.: Quantitative uniqueness properties for $L^2$ functions with fast decaying, or sparsely supported, Fourier transform. \emph{Int. Math. Res. Not. IMRN} \textbf{2022}(16),
12148--12179 (2022). \url{https://doi.org/10.1093/imrn/rnab075}.

\bibitem{JinZhang}
Jin, L. and Zhang, R.: Fractal uncertainty principle with explicit exponent. \emph{Math. Ann.} \textbf{376}(3--4), 1031--1057 (2020). \url{https://doi.org/10.1007/s00208-019-01902-8}.

\bibitem{Khats}
Khats', R.: An analogue of the Paley--Wiener theorem for the Hankel type
transform in weighted $L^2$-spaces. \emph{Gulf J. Math.} \textbf{17}(1), 1--11 (2024). \url{https://doi.org/10.56947/gjom.v17i1.1663}.

\bibitem{Koosis}
Koosis, P.: \emph{The Logarithmic Integral I}. Cambridge Studies in Advanced Mathematics, vol.~12. Cambridge University Press, Cambridge (1988).

\bibitem{MashreghiNazarovHavin}
Mashreghi, J., Nazarov, F.L. and Havin, V.P.:
Beurling--Malliavin multiplier theorem: the seventh proof.
\emph{St. Petersburg Math. J.} \textbf{17}(5), 699--744 (2006).
\url{https://doi.org/10.1090/S1061-0022-06-00926-5}.

\bibitem{VasMajorant}
Vasilyev, I.: A generalization of the Beurling--Malliavin majorant theorem. \emph{Anal. PDE} \textbf{17}(9), 3355--3369 (2024). \url{https://doi.org/10.2140/apde.2024.17.3355}.

\bibitem{Vas}
Vasilyev, I.: The Beurling and Malliavin theorem in several dimensions. \emph{Math. Ann.} \textbf{391}, 6057--6072 (2025). \url{https://doi.org/10.1007/s00208-024-03071-9}.

\bibitem{ZWD}
Zhao, X., Wei, L. and Duan, Z.: Fractal uncertainty and quantitative uniqueness for the Fourier Bessel transform. Preprint, arXiv:2608.15126 (2026). \url{https://arxiv.org/abs/2608.15126}.

\end{thebibliography}
\end{document}